\documentstyle[fullpage,12pt]{amsart}
\renewcommand{\marginpar}[1]{}
\catcode`\@=12

\def\Empty{}
\newcommand\oplabel[1]{
  \def\OpArg{#1} \ifx \OpArg\Empty {} \else
  	\label{#1}
  \fi}
		
\long\def\realfig#1#2#3#4{
\begin{figure}[htbp]
\centerline{\psfig{figure=#2,width=#4}}
\caption[#1]{#3}
\oplabel{#1}
\end{figure}}

\newcommand{\comm}[1]{}
    
\input{psfig}
\newtheorem{thm}{Theorem}[section]

\newtheorem{lem}[thm]{Lemma}
\newtheorem{prop}[thm]{Proposition}
\newtheorem{schw}{Schwarz Lemma}

\theoremstyle{definition}
\newtheorem{defn}{Definition}[section]

\theoremstyle{remark}
\newtheorem{rem}{Remark}[section]

\newcommand{\diam}{\operatorname{diam}}
\newcommand{\dist}{\operatorname{dist}}
\newcommand{\meas}{\operatorname{meas}}
\newcommand{\cl}{\operatorname{cl}}
\newcommand{\tl}{\tilde}
\newcommand{\wtl}{\widetilde}
\newcommand{\eps}{\epsilon}

\newcommand{\EE}{{\cal E}}
\renewcommand{\frak}[1]{\EE}

\numberwithin{equation}{section}
\newsymbol\Subset 1362
\newcommand{\thmref}[1]{Theorem~\ref{#1}}

\newcommand{\secref}[1]{\S\ref{#1}}
\newcommand{\lemref}[1]{Lemma~\ref{#1}}
 
\newcommand{\figref}[1]{Fig.~\ref{#1}}
\newcommand{\ang}[2]{\widehat{(#1,#2)}}
\newcommand{\C}[1]{{\Bbb C_{#1}}}

\begin{document}

\title[Complex Bounds]{Complex Bounds for Critical Circle Maps}
\author{Michael Yampolsky}
\maketitle

\def\IMSmarkvadjust{0 pt}
\def\IMSmarkhadjust{0 pt}
\def\IMSmarkhpadding{0 pt}
\def\IMSpubltext{Published in modified form:}
\def\SBIMSMark#1#2#3{
 \font\SBF=cmss10 at 10 true pt
 \font\SBI=cmssi10 at 10 true pt
 \setbox0=\hbox{\SBF \hbox to \IMSmarkhpadding{\relax}
                Stony Brook IMS Preprint \##1}
 \setbox2=\hbox to \wd0{\hfil \SBI #2}
 \setbox4=\hbox to \wd0{\hfil \SBI #3}
 \setbox6=\hbox to \wd0{\hss
             \vbox{\hsize=\wd0 \parskip=0pt \baselineskip=10 true pt
                   \copy0 \break%
                   \copy2 \break%
                   \copy4 \break}}
 \dimen0=\ht6   \advance\dimen0 by \vsize \advance\dimen0 by 8 true pt
                \advance\dimen0 by -\pagetotal
	        \advance\dimen0 by \IMSmarkvadjust
 \dimen2=\hsize \advance\dimen2 by .25 true in
	        \advance\dimen2 by \IMSmarkhadjust

%
%
  \openin2=publishd.tex
  \ifeof2\setbox0=\hbox to 0pt{}
  \else 
     \setbox0=\hbox to 3.1 true in{
                \vbox to \ht6{\hsize=3 true in \parskip=0pt  \noindent  
                {\SBI \IMSpubltext}\hfil\break
                {\it  Ergodic Th. \& Dynam. Sys.}~{\bf 19} (1999), 227--257.
 
                \vfill}}
  \fi
  \closein2
  \ht0=0pt \dp0=0pt
 \ht6=0pt \dp6=0pt
 \setbox8=\vbox to \dimen0{\vfill \hbox to \dimen2{\copy0 \hss \copy6}}
 \ht8=0pt \dp8=0pt \wd8=0pt
 \copy8
 \message{*** Stony Brook IMS Preprint #1, #2. #3 ***}
}

\SBIMSMark{1995/12}{October 1995}{}
\thispagestyle{empty}

\begin{abstract}
We use the methods developed with M. Lyubich for proving complex bounds
for real quadratics to extend E. De Faria's complex a priori bounds
to all critical circle maps with an irrational rotation number.
The contracting property for renormalizations of critical circle maps follows.
\\
In the Appendix we give an application of the complex bounds for
proving local connectivity of some Julia sets.
\end{abstract}

\section{Introduction}
The object of our consideration will be 
the family of analytic orientation-preserving 
self-homeomorphisms of the circle $f:{\bf T\to \bf T}$ with 
one critical point  at $0$. To fix our  ideas we will assume
that the  critical point has order three,
although our considerations are valid in the  general case.
 
 We will further refer
to such maps as {\it critical circle maps}.
A critical circle map has a well defined rotation number denoted
further by $\rho(f)$.
Examples of critical circle maps with any given rotation number
are provided by the maps of the standard family 
$x\mapsto x+\theta-\frac{1}{2\pi}\sin(2\pi x)$.

We will be considering only  mappings with irrational
rotation number. 
For such mappings Lanford \cite{La1,La2} has defined the infinite sequence
of renormalizations.

De~Faria in \cite{dF1,dF2} has  extended
these renormalizations to the complex plane, and has 
developed the renormalization theory for
critical circle maps, parallel to the Sullivan  theory
for unimodal maps of the interval. The key condition for applicability
of this theory is the existence of certain geometric bounds
on the renormalized maps, {\it the complex a priori bounds}.

Applying Sullivan's methods to circle mappings de~Faria has shown
the existence of complex a priori bounds for the class of mappings 
with $\rho(f)\in \operatorname{dioph}^2$.

In this paper we use the method developed with M. Lyubich \cite{LY}
to show the existence of complex a priori bounds for all critical
circle mappings with irrational $\rho$.

The method was originally used in \cite{LY} to treat a special case 
of combinatorics of quadratic maps of an interval, {\it the essentially
bounded combinatorics}. However, this particular 
case for the interval maps corresponds to the most general case 
for circle maps, which allows us to prove the following.

\begin{thm}
\label{bounds}
Let $f$ be a critical circle map of Epstein class. Then
$f$ has complex a priori bounds.
\end{thm}

The theorem relies on the following key cubic estimate for the 
renormalizations of the map  $f$:
\begin{equation}
\label{key estimate}
|R^kf(z)|\geq c|z|^3,
\end{equation}
with an absolute $c>0$. 
The proof of the key estimate is outlined in \secref{outline}.

The above result allows one to extend the renormalization theory
developed by de~Faria to all critical circle mappings (in preparation).

In the appendix we give another application of the complex a priori bounds,
a new proof of the result of C. Petersen on local connectivity of Julia sets
of some Siegel quadratics.\\ \\

{\bf\noindent Acknowledgments.} I would like to thank my academic advisor
Mikhail Lyubich for a lot of helpful conversations and suggestions.
I am grateful to Carsten Petersen for discussing some of his results with me.
I also thank Scott Sutherland for his help with computer pictures.
Finally most of this work was done at MSRI at Berkeley, and I would like to 
thank them for their hospitality.

\section{Preliminaries}
\label{prelim}
\subsection{General notations and terminology} 

We use $|J|$ for the length of an interval $J$, 
 $\dist$ and $\diam$ for the Euclidean distance and diameter in ${\Bbb C}$.
The notation $[a,b]$ stands for  the (closed) interval with endpoints $a$ and $b$ without
specifying their order. 

Two sets  $X$ in $Y$ in ${\Bbb C}$ are called {\it $K$-commensurable}
or simply {\it commensurable} if 
$$K^{-1}\leq \diam X/\diam Y\leq K$$ with a constant $K>0$ which 
may depend only on the specified combinatorial bounds.
A configuration of points $x_1,\dots\,x_n$ is called {\it $K$-bounded}
if any two intervals $[x_i,x_j]$, and $[x_k,x_l]$ are $K$-commensurable.

We say that an annulus $A$ has a {\it definite modulus}
 if $\mod A\geq \delta>0$,
where $\delta $  depends only on the specified  bounds.  

For a pair of intervals $I\subset J$ we say that $I$ is contained 
{\it well inside} of $J$ if for each component $L $ of $J\setminus I$
we have
 $|L|\geq K|I|$ where the constant $K>0$ depends only on the specified
bounds.

Let $f$ be an  analytic map  whose restriction to the
circle $T$ is a critical circle map.  We reserve
the notation $f^{-i}(z)$ for the $i-$th preimage of $z\in T$ under 
$f|_T$.

\subsection{Hyperbolic disks}
Given an interval $ J\subset {\Bbb R}$, let
$\C{J}\equiv {\Bbb C}\backslash ({\Bbb R}\backslash {\rm J})$ denote
the plane slit along two rays. Let  $\overline{{ {\Bbb C}}_J}$
 denote the completion of
this domain in the path metric in ${\Bbb C}_{J}$ (which means that we add to
 ${\Bbb C}_{J}$
 the banks of the slits).

By symmetry, $J$ is a hyperbolic geodesic in ${\Bbb C}_{J}$.
The {\it geodesic neighborhood} of $J$  of radius
$r$ is the 
set of all points in ${\Bbb C}_{J}$
 whose hyperbolic distance to $J$ 
 is less than $r$. It is easy to see that such a neighborhood is
the union of two $\Bbb R$-symmetric segments of Euclidean disks 
based on $J$
 and having angle $\theta=\theta(r)$ with $\Bbb R$. Such a hyperbolic 
disk will be denoted by $D_{\theta}(J)$ (see Figure 1).
Note, in particular, that the Euclidean disk $D(J)\equiv D_{\pi/2}(J)$
can also be interpreted as a hyperbolic disk.

\realfig{gdang}{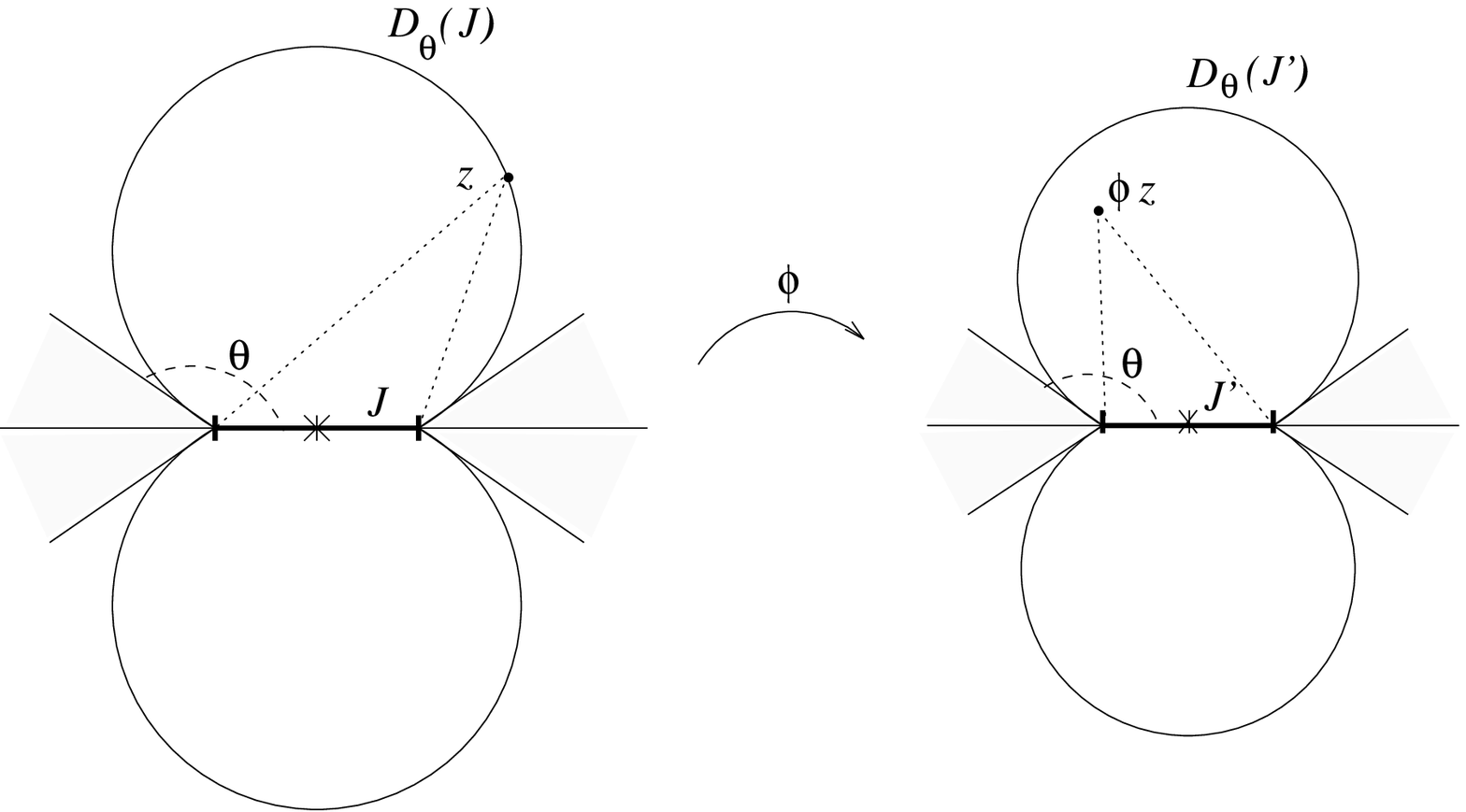}{}{15cm}

These hyperbolic neighborhoods were introduced into the subject by Sullivan
[S]. They are a key tool for getting complex bounds due to
the following version of the Schwarz Lemma:

\begin{schw}
  Let us consider two intervals  
$J'\subset J\subset {\Bbb R}$. 
Let $\phi:{\Bbb C}_{J}\rightarrow{\Bbb C}_{J'}$ be an analytic map such that
$\phi(J)\subset J'$. Then for any $\theta\in (0, \pi),$ 
$\phi (D_\theta (J))\subset D_\theta (J')$.
\end{schw}

Let $J=[a,b]$.  
For a point $z\in {\bar {\Bbb C}}_J$, the 
{\it angle between $z$ and 
 $J$}, $\ang{z}{J}$
is the least of the angles  between the intervals $[a,z]$, $[ b, z]$  and 
the corresponding rays $(a,-\infty]$, $[b, +\infty)$ of the real line,
measured in the range $ 0\leq \theta\leq \pi$. 

 We will use the following observation to control the expansion of the inverse
branches.

\begin{lem}\label{good angle}
Under the assumptions of the Schwarz Lemma, let us consider 
 a point $z\in {\Bbb C}_{J}$ such that $\dist(z, J)\geq |J|$ and 
$\ang{z}{J}\geq\epsilon$. Then
$${\dist(\phi z, J')\over |J'|}\leq C{\dist(z, J)\over |J|}$$ 
for some constant $C=C(\epsilon)$
\end{lem}

\begin{pf} Let us normalize the situation in this way: $J=J'=[0,1]$.
Notice that the smallest (closed) geodesic neighborhood 
$\cl D_{\theta}(J)$ enclosing $z$ satisfies:\\ 
$\diam D_{\theta}(J)\leq C(\epsilon)\dist (z, J)$ (cf \figref{gdang} ).

Indeed, if $\theta\geq \epsilon/2$ then 
$\diam D_{\theta} (J)\leq C(\epsilon)$, which is fine since
$\dist (z, J)\geq 1$.\\
Otherwise the intervals $[0,z]$ and $[1,z]$ cut out sectors of angle 
size at least $\epsilon/2$ on the circle $\partial D_{\theta}(J)$.
Hence the lengths of these intervals are commensurable with 
$\diam D_{\theta}(J)$ (with a constant depending on $\epsilon$).
Also, by elementary trigonometry these lengths are at least 
$\sqrt{2}\dist(z, J)$, provided that $\dist (z, J)\geq |J|$.

By Schwarz Lemma, 
 $\dist(\phi z, J')\leq\diam( D_{\theta}(J')),$ and the claim follows.
\end{pf}
\realfig{fig11}{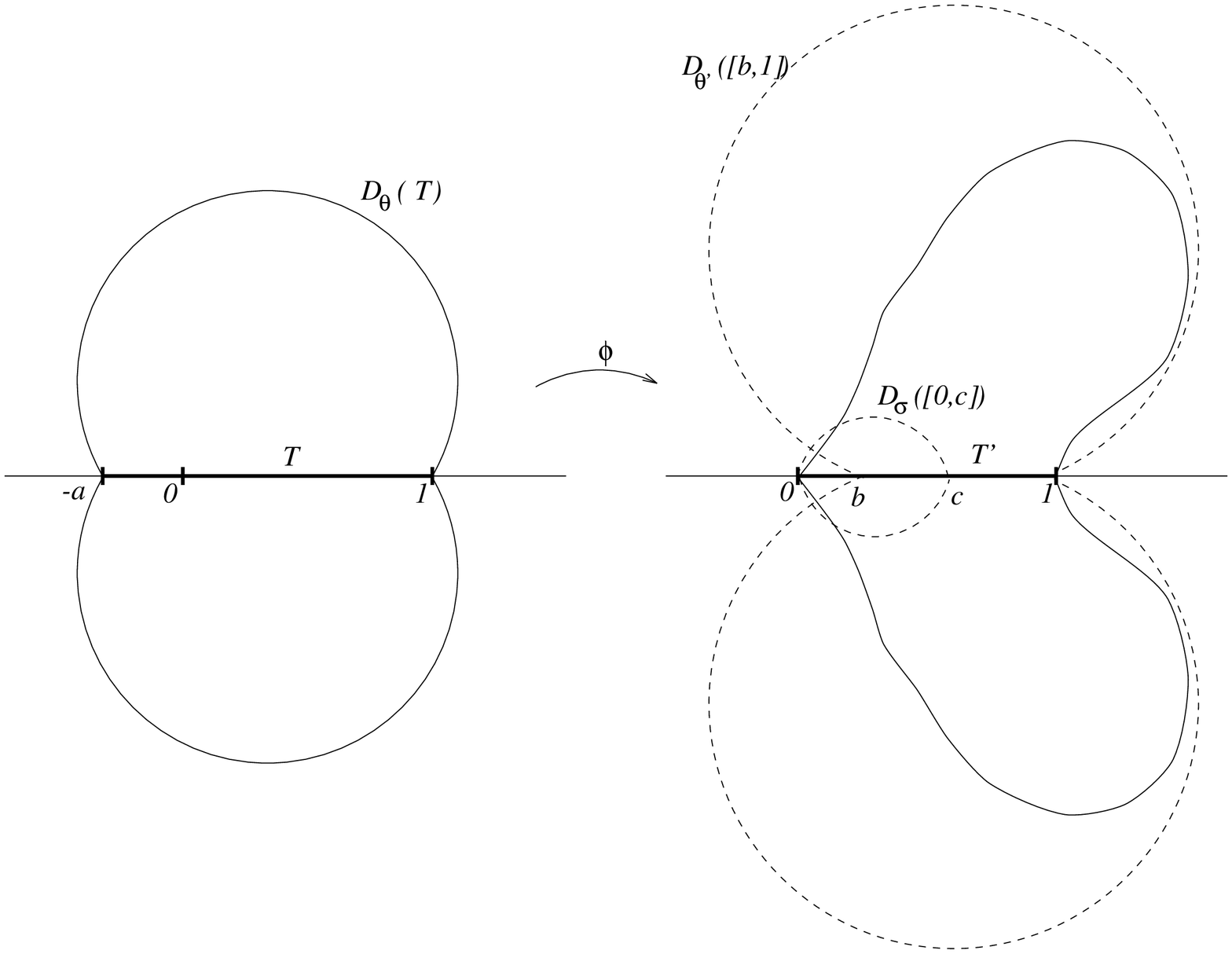}{}{14cm}

\subsection{Cube root}
In the next lemma we collect for future reference some elementary properties of the
cube root map. 
Let $\phi(z)$ be the branch of the cube root mapping the slit plane
${\Bbb C}\setminus {\Bbb R}_{-}$ into $\{ z | |\operatorname{arg}(z)|<\pi/3\}.$

\begin{lem}
\label{cube root}
let $K>1,\; \delta>0, \; K^{-1}\leq a\leq K,\; T=[-a,1],\; T'=[0,1].$
Then:
\itemize{
\item $\phi D_\theta(T)\subset D_{\theta'}(T'),$ with $\theta'$ depending on $\theta$ and $K$ only.
\item Moreover, there exist $b,c\in[0,1]$, such that $0,b,c,1$ form a 
$C(K)$-bounded configuration, and $\phi D_\theta(T)\subset 
D_{\theta'}([b,1])\cup D_\sigma([0,c])$ for $\sigma<\pi/2$ cf. \figref{fig11}.
\item If $z'\in \phi D(T)\setminus D([-\delta,1+\delta])$, then
$$ \ang{z'}{T'}>\eps(K,\delta)>0\; {\rm and}\; C(K,\delta)^{-1}<\dist(z',T')<C(K,\delta).$$
}
\end{lem}

%
%
\comm{

\subsection{ Renormalization.}
\label{renorm}
Denote by $q_n$ the moments of closest returns of the critical point $0$.
Let $I_m\equiv [0,f^{q_m}(0)]$.
By \'Swiatek- Herman real a priori bounds (\cite{Sw,H}), the
intervals $I_m$ and $I_{m+1}$ are $K$-commensurable with a universal
constant $K$ provided $m$ is large enough.

The pair of maps $({\tl f}^{q_{m+1}}|{\tl I}_m;\; {\tl f}^{q_m}|{\tl I}_{m+1})$
where $\tl{}$ means rescaling by a linear factor $\lambda=\frac
{1}{|f^{q_m}(0)|}$, forms a (rescaled) {\it real commuting pair}.
Its {\it renormalization} is the rescaled pair
$({\tl f}^{q_{m+2}}|{\tl I}_{m+1};\; {\tl f}^{q_{m+1}}|{\tl I}_{m+2})$.
}
\subsection{Real Commuting Pairs and Renormalization}
\label{renorm}

We present here a brief summary on renormalization of critical circle
mappings. A more extensive exposition can be found for example
in {\cite{dF2}.

Let $f$ be a critical circle mapping with an irrational rotation number
$\rho(f)$. Let $\rho(f)$ have an infinite continued fraction expansion 
$$\rho(f)=[r_0,r_1,r_2,\ldots]=\cfrac{1}{r_0+\cfrac{1}{r_1+\cfrac{1}{r_2+\dotsb}}}.$$

We say that $\rho$ is of {\it bounded type} if $\sup r_i <\infty$.
This is equivalent to $\rho\in \operatorname{dioph}^2$.

Denote by $q_m$ the moments of closest returns of the critical point 
$0$. Note that the numbers $q_m$ appear as the denominators in the
irreducible form of the $m$-th truncated continued fraction expansion
of $\rho(f)$, ${p_m\over q_m}=[r_0,\ldots,r_{m-1}]$.

Let $I_m\equiv[0,f^{q_m}(0)]$. By \'Swiatek- Herman real a priori bounds
(\cite{Sw,H}), the intervals $I_m$ and $I_{m+1}$ are $K-$commensurable, 
with a universal constant $K$ provided $m$ is large enough.

The dynamical first return map on the interval $I_m\cup I_{m+1}$ is
$f^{q_m}$ on $I_{m+1}$  and $f^{q_{m+1}}$ on $I_m$.
The consideration of pairs of maps 
\begin{center}
$(f^{q_{m+1}}|I_m,f^{q_m}|I_{m+1})$
\end{center}
leads to the following general definition due to
Lanford and Rand (cf. \cite{La1,La2,Ra1,Ra2}).
 
\begin{defn}
A (real) commuting pair $\zeta=(\eta,\xi)$ consists of two real
orientation preserving smooth homeomorphisms 
$\eta:I_\eta\to{\Bbb R},\;\xi:I_{\xi}\to{\Bbb R}$, where
\begin{itemize}
\item{$I_\eta=[0,\xi(0)],\; I_\xi=[\eta(0),0]$;}
\item{Both $\eta$ and $\xi$ have homeomorphic extensions to interval
neighborhoods of their respective domains which commute, i.e.
$\eta\circ\xi=\xi\circ\eta$ where both sides are defined;}
\item{$\xi\circ\eta(0)\in I_\eta$;}
\item{$\eta'(x)\ne 0\ne \xi'(y) $, for all $x\in I_\eta\setminus\{0\}$,
 and all $y\in I_\xi\setminus\{0\}$.}
\end{itemize}
\end{defn}

A {\it critical commuting pair} is a commuting pair $(\eta,\xi)$,
which maps can be decomposed as $\eta=h_\eta\circ Q\circ H_\eta$,
and $\xi=h_\xi\circ Q\circ H_\xi$, where $h_\eta,h_\xi,
H_\eta,H_\xi$ are real analytic  diffeomorphisms and $Q(x)=x^3$.

Given a commuting pair $\zeta=(\eta,\xi)$ we will denote by 
$\wtl\zeta$ the pair $(\wtl\eta|\wtl{I_\eta},\wtl\xi|\wtl{I_\xi})$
where tilde  means rescaling by a linear factor $\lambda={1\over |I_\eta|}$.


For a critical circle mapping $f$ one obtains a critical commuting pair
from the pair of maps $(f^{q_{m+1}}|I_m,f^{q_m}|I_{m+1})$ as follows.
Let $\bar f$ be the lift of $f$ to the real line satisfying $\bar f '(0)=0$,
and $0<\bar f (0)<1$. For each $m>0$ let $\bar I_m\subset {\Bbb R}$ 
denote the closed 
interval adjacent to zero which projects down to the interval $I_m$.
Let $\tau :{\Bbb R}\to {\Bbb R}$ denote the translation $x\mapsto x+1$.
Let $\eta :\bar I_m\to {\Bbb R}$, $\xi:\bar I_{m+1}\to {\Bbb R}$ be given by
$\eta\equiv \tau^{-p_{m+1}}\circ\bar f^{q_{m+1}}$,
$\xi\equiv \tau^{-p_m}\circ\bar f^{q_m}$. Then the pair of maps
$(\eta|\bar I_m,\xi|\bar I_{m+1})$ forms a critical commuting pair
corresponding to $(f^{q_{m+1}}|I_m,f^{q_m}|I_{m+1})$.
Henceforth, we shall abuse notation and use
\begin{equation}
\label{real1}
(f^{q_{m+1}}|I_m,f^{q_m}|I_{m+1})
\end{equation}
to denote this commuting pair.

We see  that  return maps for critical 
circle homeomorphisms give rise to a sequence of critical commuting
pairs (\ref{real1}).
 Conversely, regarding $I_\eta$ as a circle (identifying $\xi(0)$
and $0$), we can recover a smooth conjugacy class of critical
circle mappings $f^\phi=\phi\circ f_\zeta\circ\phi^{-1}$,
where $\phi:I_\eta\to I_\eta$ is a smooth orientation preserving
homeomorphism with a cubic critical point at $0$, and
 $f_\zeta$ is a circle homeomorphism defined by

\begin{equation}
f_\zeta(x)=\left\{
\begin{array}{ll}
\xi\circ\eta(x),& {\rm if}\; 0\leq x\leq\eta^{-1}(0)\\
\eta(x),& {\rm if}\; \eta^{-1}(0)\leq x\leq \xi(0)
\end{array}\right.
\end{equation}

Let $\rho(\zeta)=\rho(f_\zeta)$ be the {\it rotation number of
the commuting pair $\zeta$}.
If $\rho(\zeta)=[r,r_1,r_2,\ldots]$, one verifies that the 
mappings $\eta|[0,\eta^r(\xi(0))]$ and $\eta^r\circ\xi|I_\xi$
again form a commuting pair.

\begin{defn}
The {\it renormalization} of a real commuting pair $\zeta=(\eta,
\xi)$ is the commuting pair
\begin{center}
${\cal{R}}\zeta=(
\widetilde{\eta^r\circ\xi}|
 \widetilde{I_{\xi}},\; \widetilde\eta|\widetilde{[0,\eta^r(\xi(0))]}).$
\end{center}
\end{defn}

It is easy to see that renormalization acts as a Gauss map on rotation
numbers, i.e. if $\rho(\zeta)=[r,r_1,r_2,\dots]$, $\rho({\cal{R}}\zeta)=
[r_1,r_2,\ldots]$.

Finally note, that the renormalization of the real commuting pair
(\ref{real1}) is the rescaled pair 
$(\wtl{{f}^{q_{m+2}}}|\wtl{{I}_{m+1}},\wtl{{f}^{q_{m+1}}}|\wtl{{I}_{m+2}})$.

In this way for a given critical circle mapping with an irrational rotation
number we obtain an infinite sequence of renormalizations
$\{ (\wtl{f^{q_{i+1}}}|\wtl{I_i},\wtl{f^{q_{i}}}|\wtl{I_{i+1}})\}_{i=1}^{\infty}$.

\subsection{Holomorphic Commuting Pairs}

\realfig{circ3}{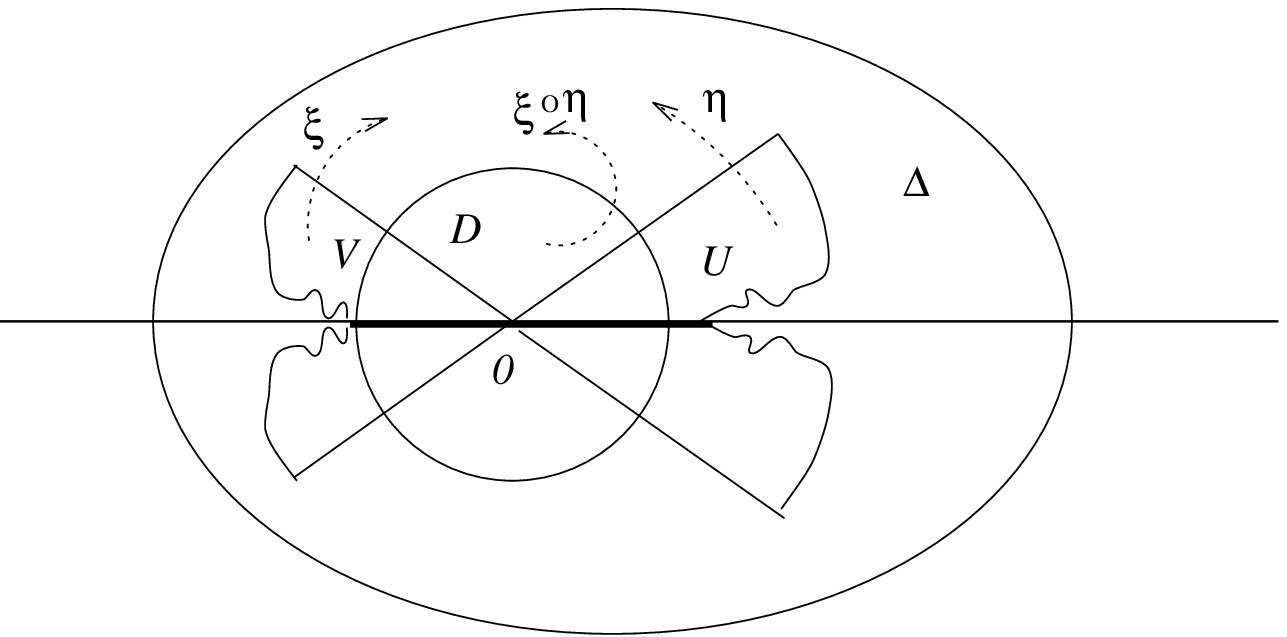}{}{12cm}
Following \cite{dF1,dF2} we say that a real commuting pair
$ ( \eta , \xi )$
extends to a {\it holomorphic commuting pair} (cf. \figref{circ3})
if there exist four 
$\Bbb R$-symmetric domains $\Delta$, $D$, $U$, $V$, such that 
\begin{itemize}
\item  $\bar D,\; \bar U,\; \bar V\subset \Delta$,
 $\bar U\cap \bar V=\{ 0\}$;
 $U\setminus D$,  $V\setminus D$, $D\setminus U$, and $D\setminus V$ 
 are nonempty 
connected sets, $U\supset I_\eta$, $V\supset I_{\xi}$;
\item mappings $\eta:U\to \Delta\cap {\Bbb C}_{\eta(J_{U})}$ and
 $\xi:V\to \Delta\cap{ \Bbb C}_{{\xi}(J_{V})}$ are onto and
univalent, where $J_{U}=U\cap {\Bbb R}$, $J_{V}=V\cap {\Bbb R}$;
\item $\eta$ and $\xi$ have holomorphic extensions to $D$ which
commute, $\eta\circ\xi(z)=\xi\circ\eta(z)$  $\forall z\in D$; 
$\eta\circ\xi:D\to \Delta\cap{\Bbb C}_{{\eta\circ\xi}(J_{D})}$,
where $J_{D}=D\cap {\Bbb R}$,
is a three-fold branched covering with the only critical point at $0$
.
\end{itemize}
Note that $J_{D}=(\eta^{-1}(0),\xi^{-1}(0))$.

\begin{defn}[Complex Bounds]
\label{cab}
We say that a critical circle map $f$ has {\it complex a priori bounds} if
there exists $M$ and $\mu >0$, such that for all $m>M$
the real commuting pair (\ref{real1}) extends to a holomorphic commuting pair
$(\Delta_m, D_m, U_m, V_m)$ with $\mod (\Delta_m\setminus (D_m\cup U_m
\cup V_m))>\mu>0$.
\end{defn}


\subsection{Epstein class}

A map $g|I$ belongs to an
 {\it Epstein class},
 if the restriction of $g$ to $I$ can be decomposed as
 $g\equiv h\circ Q$, where  
$Q:I\to J\equiv g(I)$ is a real cubic polynomial,
 and $h:J\to J$
is an orientation preserving diffeomorphism, which inverse $h^{-1}$
extends to a univalent mapping ${\Bbb C}_{\tl J}\to {\Bbb C}$,
where $\tl J\supset J$.
A commuting pair $(\eta,\xi)$
belongs to Epstein class if both maps do.

We will denote by $ {\cal E}_s$  the family of maps of Epstein class 
for which the length of each component of $\tl J\setminus J$ is
$s|J|$.

One immediately obtains:
\begin{lem} 
An Epstein class ${\cal E}_s$ is invariant under renormalization.
\end{lem}

We supply the space of maps of Epstein class with {\it Caratheodory 
topology} (see \cite{McM1}).

The next Lemma follows from a standard normality argument.

\begin{lem}
For each $s>0$ the space ${\cal E}_s$ is compact.
\end{lem}

By real a priori bounds
 there exists a universal $s>0$ such that
$(f^{q_{m+1}}|I_m;f^{q_m}|I_{m+1})\in {\cal E}_s$, for all sufficiently
large $m$, provided $f$ belongs to an Epstein class.

{\it All maps in this paper are assumed to belong to some Epstein class 
${\cal E}_s$}. 
The bounds $M$ and $\mu $ in the Definition \ref{cab} will depend only on
the chosen value of $s$. Note that even for $\rho(f)\in \operatorname{dioph}^2$
this improves the complex a priori bounds of E. de Faria which depend also on 
the combinatorics of the map.

\section{Outline of the proof}
\label{outline}

\subsection{Main Lemma}
Fix $n$, and let $p=q_{n+1}$.


Let us consider the decomposition
\begin{equation}
\label{decomposition}
f^p=\psi_n\circ f,
\end{equation}

where $\psi_n$ is a univalent map from a neighborhood of $f(I_n)$ onto
${\Bbb C}_{f^p(I_n)}$.

\begin{lem}
\label{lin}

 There exist a disc $D_0$ around $0$ and universal constants $C_1$ 
and $C_2$ depending only 
on real a priori 
bounds, such that 
$\forall z\in{\Bbb C}_{f^p(I_n)}\cap D_0$
 the following estimate holds:

\begin{equation}
\label{linear}  
\frac{\dist(\psi_n^{-1}z,f(I_n))}{|f(I_n)|}\leq C_1\left(\frac{\dist(z,I_n)}
{|I_n|}\right) +C_2,
\end{equation}

where $\psi_n$ is the map from (\ref{decomposition}).
\end{lem}

Thus the map $\psi^{-1}_n$ has at most linear growth.

Note that if $\ang{z}{I_n}>\eps>0$, then the inequality (\ref{linear})
follows directly from \lemref{good angle} with the constants depending only on
$\eps$. Our strategy of proving \lemref{lin} will be to monitor the 
inverse orbit of a point $z$ together with the interval $I_n$ until
they satisfy this ``good angle'' condition.

\lemref{lin} immediately yields the key cubic estimate:
\begin{prop}[Key estimate]
\label{cubic-estimate}
Let $U$ be the neighborhood of $\widetilde{ I_n}$ which is mapped univalently 
by $\widetilde{f^p}$
onto $\C{\widetilde{ f^{p}}(\widetilde{ I_n})}$. There exist universal 
constants $B$ and $C$, and a disc $D_0$ around $0$ 
such that $\forall z\in U$ with   $\wtl{f^p}(z)\in \widetilde{D_0}$ and
 $|\wtl{f^p}(z)|>B$, one has:
\begin{center}
$|\widetilde{ f^p}(z)|>C|z|^3$.
\end{center}
\end{prop}

The \thmref{bounds} now follows. 

\begin{rem}
\label{commens}
We remark that it follows from our estimate that the domains
 $\Delta_m, D_m, U_m,$ and $V_m$ in the definition
of complex bounds (\ref{cab}) can be chosen $K$-commensurable with
$I_m$, with a universal constant $K$.
\end{rem}

\section{Proof of \lemref{lin}}

\subsection{}
Let $g:U\to \C{T}$ be a map of Epstein class. Take an $x\in {\Bbb R}$, and
$z\in \C{T}$. If we have a backward orbit of $x\equiv x_0,x_{-1},\ldots,
x_{-l}$ of $x$ which does not contain $0$, the {\it corresponding }
backward orbit $z\equiv z_0,z_{-1},\ldots,z_{-l}$ is obtained by applying
the appropriate branches of the inverse functions: $z_{-i}=g_{x_{-i}}(z)$,
where $g_{x_{-i}}(x)=x_{-i}$.

Let $D_m$ denote the disc $D([f^{q_{m+1}}(0),f^{q_{m}-q_{m+1}}(0)])$.

Let $p$ be as above, and
 consider the inverse orbit:
\begin{equation}
\label{J-orbit}
J_0\equiv f^p(I_n),J_{-1}\equiv f^{p-1}(I_n),\dots,
J_{-(p-1)}\equiv f(I_n)
\end{equation}
For a point $z\in {\Bbb C}_{T_n}$ consider the corresponding inverse orbit
\begin{equation}
\label{z-orbit}
z_0\equiv z,z_{-1},\dots,z_{-(p-1)}
\end{equation}

We say that a point $z_{-i}$ $\eps$-{\it jumps} if $\ang{z_{-i}}{J_{-i}}>\eps$
for some fixed $\eps >0$.

Let us call the moment $-i$ ``good'' if the interval $J_{-i}$ is 
commensurable with $J_0$. For example the first few returns of the
orbit (\ref{J-orbit}) to any of the $I_m$ are good.

We would like to assert that the points of the orbit (\ref{z-orbit}) either
$\eps$-jump at a good moment, or follow closely the corresponding intervals
of the orbit (\ref{J-orbit}).

The first step towards this assertion is the following 

\begin{lem}
\label{cb1}
Let $J\equiv J_{-k},J_{-k-q_{m+1}}\equiv J'$ be two consecutive returns of the 
backward orbit (\ref{J-orbit}) to $I_m$, and let $\zeta$ and
$\zeta'$ be the corresponding points of the orbit (\ref{z-orbit}).


Suppose $\zeta\in D_m$, then 
%
either $\zeta'\in D_m$, or
$\ang{\zeta'}{J'}>\eps$, and $\dist({\zeta'},{J'})<C|I_m|$.
\end{lem}
\realfig{cir1}{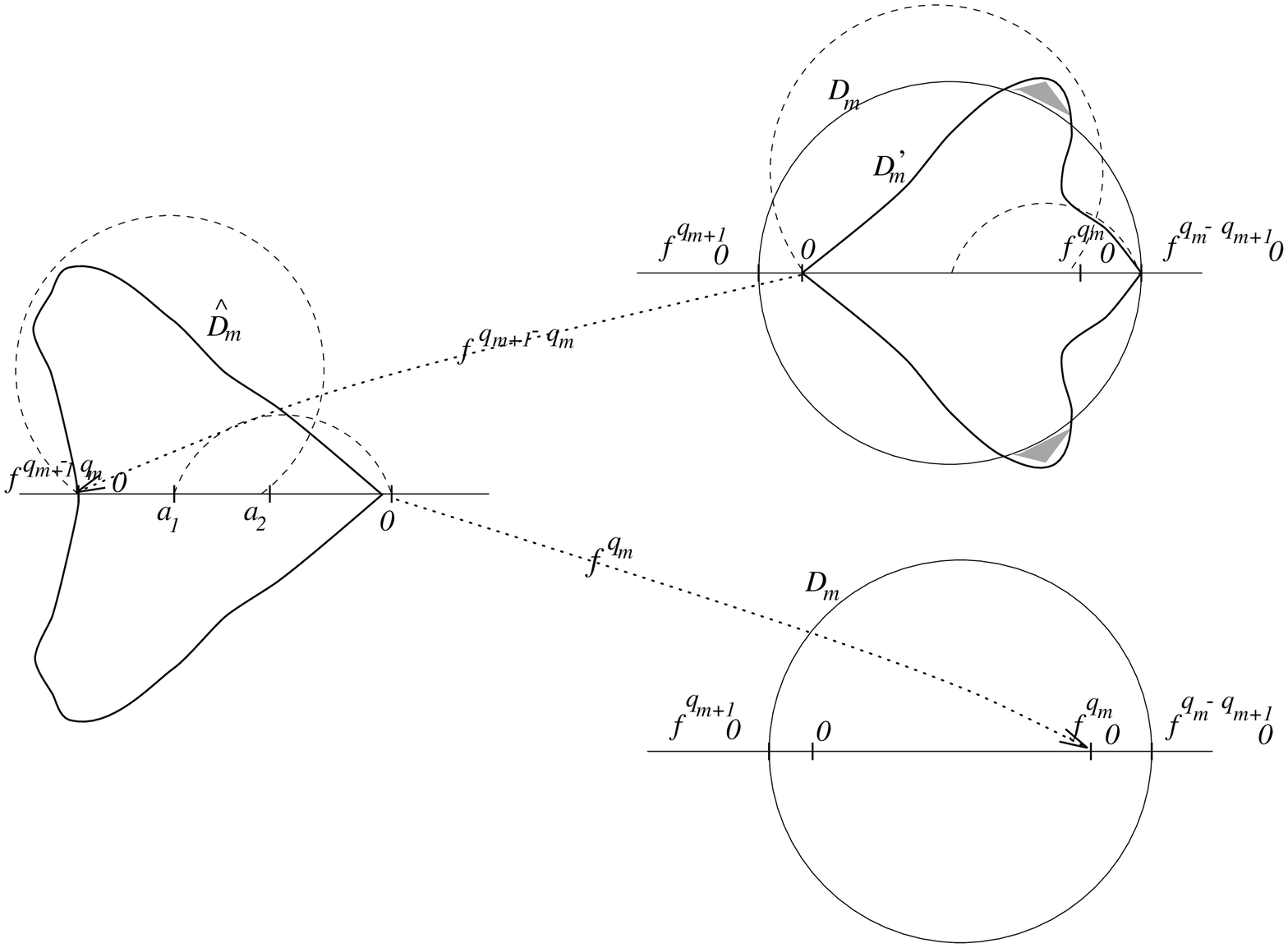}{}{17cm}
\begin{pf}
Let $D'_m$ denote the pull back of $D_m$ corresponding to the piece
of backward orbit $J_{-k},\dots,J_{-k-q_{m+1}}$, and let
$\hat{D}_m$ denote the pull back of $D_m$ along the piece of the orbit
$J_{-k}\to\dots\to J_{-k-q_m}$  (cf. \figref{cir1}).
By Schwarz Lemma and by \lemref{cube root}, there exist  points $a_1, a_2\in
[f^{q_{m+1}-q_m}(0),0]$, such that  $f^{q_{m+1}-q_m}(0), a_1, a_2$, 
and $0$ is  a $K$-bounded configuration for some
$K$ independent on $m$,  and angles
$\theta$ and $\sigma <\pi/2$, also independent on $m$, such that 
$\hat{D}_m\subset D_\theta([f^{q_{m+1}-q_m}(0),a_1])\cup D_\sigma([a_2,0])$.
Applying Schwarz Lemma, we obtain that $D'_m\subset D_m\cup D_\theta 
([0,f^{-q_{m+1}+q_{m}}(a_1)])$ and the claim immediately follows.
\end{pf}

\subsection{Saddle-node phenomenon}
Let us note first that when $q_{m+1}/q_{m}$ is large, the map $f^{q_{m+1}}:
I_m\to I_m\cup I_{m+1}$ is a small perturbation of a map with a parabolic 
fixed point ( see e.g. \cite{H} ). The next lemma 
is a direct consequence of real bounds.

\begin{lem}
\label{perturb}
Consider a family of maps 
 $\{ f_j^{q_{m_j+1}}|I_{m_j}\}$ with $f_j\in \EE_s$ and $q_{m_j+1}/q_{m_j}
\to\infty$. Then any limit point for this sequence
in the Caratheodory topology has a parabolic fixed point.
\end{lem}

We would like to deal with the ``dangerous'' situation when the backward 
orbit (\ref{z-orbit}) leaves $D_m$ at a ``bad'' moment when $J_{-i}$ is
not commensurable with the original size. The next lemma takes care of
this possibility.

\begin{lem}
\label{parab}
Let us consider the map $f^{q_{m+1}}|I_m$. Let $P_0,P_{-1},\dots,P_{-k}$
be the consecutive returns of the backward orbit (\ref{J-orbit}) to $I_m$,
and denote by $\zeta_{0},\dots,\zeta_{-k}$ the corresponding moments of
the backward orbit of a point $\zeta_0=z\in D_m$. If the ratio $q_{m+1}/q_m$
is sufficiently big, then either $z'\equiv \zeta_{-k}\in D_m$, or 
$\ang{z'}{P_{-k}}>\eps$ and $\dist(z',P_{-k})\leq C|I_m|$
\end{lem}

\begin{pf}
To be definite, let us assume that the intervals $P_{-i}$ lie on the left
 of $0$. Without loss of generality we can assume that $z\in {\Bbb H}$.
Let $\phi =f^{-q_{m+1}}$ be the branch of the inverse for which $\phi P_{-i}=
P_{-(i+1)}$. As $\phi$ is orientation preserving on $(-\infty, f^{q_{m+1}+q_{m}}(0)]$,
it maps the upper half-plane ${\Bbb H}$ into itself : 
$\phi({\Bbb H})\subset \{z=re^{i\theta}|r>0,\pi>\theta >2\pi/3\}$.

By \lemref{perturb}, if $q_{m+1}/q_m$ is sufficiently large, the map $\phi$
has an attracting fixed point $\eta_\phi\in D(I_m)\subset D_m$ (which is a 
perturbation of a parabolic fixed point). By the Denjoy-Wolf Theorem,
$\phi^n(\zeta)\underset{n\to\infty}{\to}\eta_\phi$ for any $\zeta\in {\Bbb H}$,
uniformly on compact sets of $\Bbb H$. Thus for a given compact set $K\Subset
{\Bbb H}$, there exists $N=N(K,\phi)$ such that $\phi^N(K)\subset
D_m$. By a normality argument the choice of $N$ is independent of a particular 
$\phi$ under consideration.

Suppose $\zeta_{-i}\notin D_m$. By \lemref{cb1} the set $K=(D_m\setminus
\phi(D_m))\cap {\Bbb H}$ is compactly contained in ${\Bbb H}$, and $\diam K<
C|I_m|$. For $N$ as above we have $z'\in \cup^{N-1}_{i=0}\phi^i(K)\cup D_m$,
and the lemma is proved.
\end{pf}

\subsection{The inductive step}

The next lemma provides us with the inductive steps along the 
intervals $I_m$.
\begin{lem}
\label{cb2}
Let $J$ be  the last return of the backward  orbit (\ref{J-orbit})
to the interval $I_m$ before the first return to $I_{m+1}$,
and let $J'$ and $J''$ be the first two returns of (\ref{J-orbit})
 to $I_{m+1}$. Let $\zeta$, $\zeta'$, $\zeta''$ be the corresponding moments
in the backward orbit (\ref{z-orbit}),  $\zeta =f^{q_m}(\zeta')$, 
$\zeta' =f^{q_{m+2}}(\zeta'')$.

Suppose $\zeta\in D_m$. Then either $\ang{\zeta''}{I_{m+1}}>\eps$,
or $\zeta''\in D_{m+1}$.
\end{lem}

\realfig{circ2}{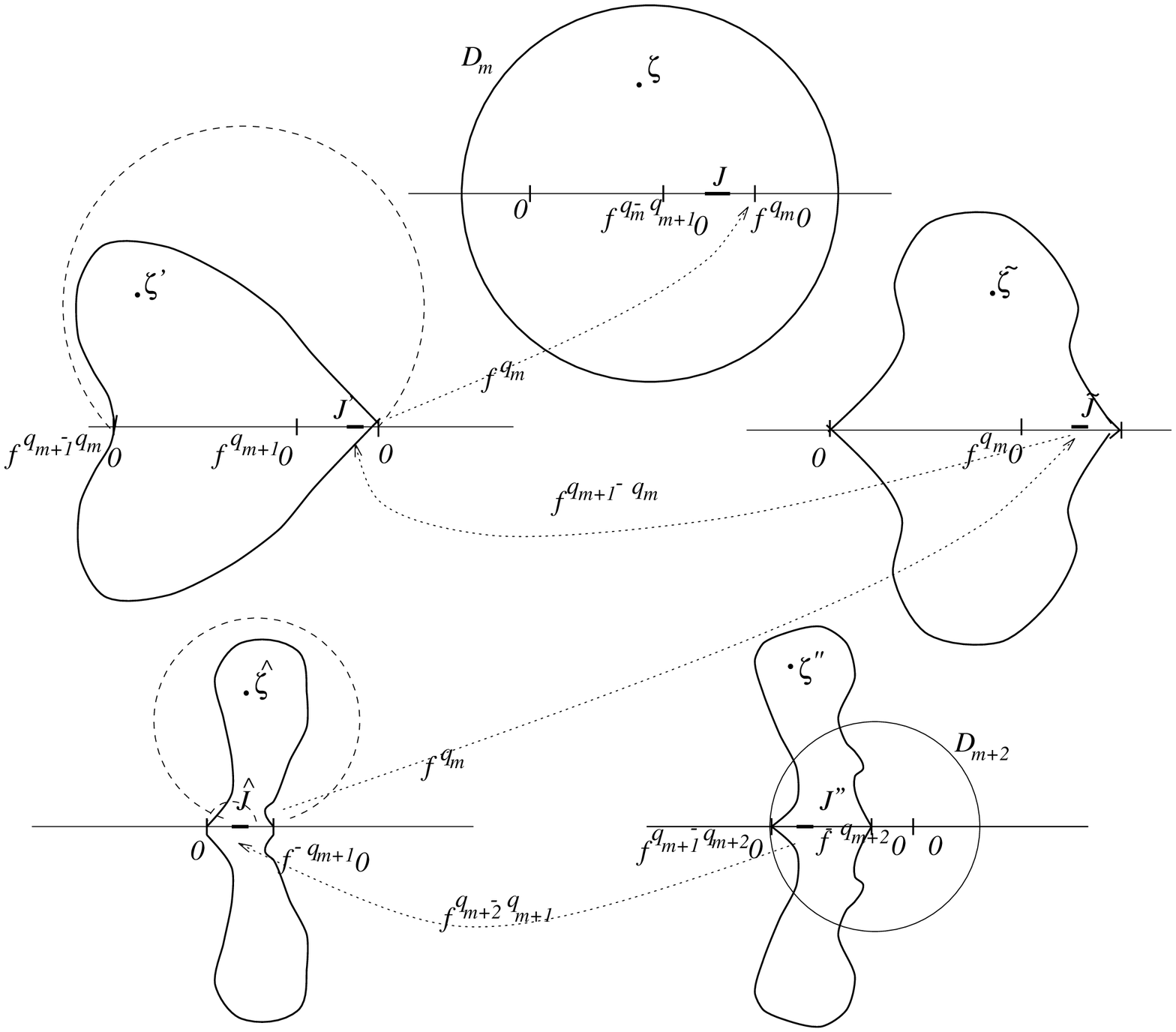}{}{17cm}
\begin{pf}
Note that $J\subset [f^{q_m+q_{m+1}}(0),f^{q_m}(0)]$.
By  \lemref{cube root}, 
 $\zeta' \in D_\theta ([f^{q_{m+1}-q_m}(0),0])$ for some 
uniform constant $\theta$. Denote by $\tl{J}$, and $\hat{J}$ the intervals
of (\ref{J-orbit}), such that $f^{q_{m+1}-q_m}(\tl J)=J'$, 
$f^{q_m}(\hat J)=\tl J$, and let $\tl \zeta$, $\hat \zeta$ be the corresponding
points in the orbit (\ref{z-orbit}) (cf. \figref{circ2}) .
The interval $\tl J\subset [f^{q_m}(0),f^{q_{m}-q_{m+1}}(0)]$, $\tl{\zeta} \in 
D_\theta ([f^{q_m}(0),f^{q_{m}-q_{m+1}}(0)])$.
 By Schwarz Lemma and \lemref{cube root} there are points $b_1,b_2\in [0,
f^{-q_{m+1}}(0)]$, such that $0, b_1, b_2,$ and  $f^{-q_{m+1}}(0)$
form  a $K$-bounded configuration, and $\hat {\zeta}\in D_\sigma ([0,b_1])\cup
D_\gamma ([b_2,f^{-q_{m+1}}(0)])$ for uniform $\gamma$, and $\sigma<\pi/2$.
The claim now follows by Schwarz Lemma.
\end{pf}

\subsection{Inductive argument}
\label{ia}

We start with a point $z\in D_1$. Consider the largest $m$ such that
$D_m$ contains $z$. We will carry out induction in $m$.
Let $P_0,\dots,P_{-k}$ be the consecutive returns of the backward orbit
(\ref{J-orbit}) to the interval $I_m$ until the first return to $I_{m+1}$,
and denote by $z=\zeta_0,\dots,\zeta_{-k}=\zeta'$ the corresponding points
of the orbit (\ref{z-orbit}).
By \lemref{cb1} and \lemref{parab}, $\zeta_{-i}$ either $\eps$-jumps
at a good moment when $P_{-i}$ is commensurable with $J_0$, and 
$\dist(\zeta_{-i},P_{-i})\leq C|I_m|$, or 
$\zeta'\in D_m$.

 In the former case we are done by \lemref{good angle}.
In the latter case consider the point $\zeta''$ which corresponds to the
second return of the orbit (\ref{J-orbit}) to $I_{m+1}$. By \lemref{cb2},
either $\ang{\zeta''}{I_{m+1}}>\eps$, and $\dist({\zeta''},{I_{m+1}})\leq
C|I_{m+1}|$, or $\zeta''\in D_{m+1}$.

In the first case we are done again by \lemref{good angle}. In the 
second case, the argument is completed by induction in $m$.

The following Remark is used in the Appendix.

\begin{rem}
\label{choice-domain}
By the above argument, in \lemref{lin} and Proposition \ref{cubic-estimate},
we can let $D_0=D_1$. We can choose instead $D_0=D_\alpha([f^{q_2}(0),
f^{q_1-q_2}(0)])$ for some $\alpha>\frac{\pi}{2}$, after an obvious change
in the argument.
\end{rem}

\appendix
\section{Application to Proving Local Connectivity}
   
\subsection{Preliminaries}
 
Let $\Bbb D$ denote the unit disc, let $T=\partial {\Bbb D}$.

For a measurable set $S\subset {\Bbb C}$ let $\operatorname{meas}(S)$
denote its planar Lebesgue measure.

\comm{
For an irrational number $\theta$ consider its continued fraction expansion
\begin{center}
$ \theta = [a_0,a_1,a_2,\ldots,a_i,\ldots]$.
\end{center}
 We say that $\theta$ is of
{\it constant type} if $\sup a_i<\infty$. Note that this is equivalent
to $\theta\in \operatorname{dioph}^2$.
}

Let $P_\theta(z)=e^{i2\pi\theta}z+z^2$. Denote by $f_\theta$ the 
Blaschke product 
\begin{equation}
\label{blyashke}
f_\theta(z)=e^{i2\pi\tau(\theta)}z^2\frac{z-3}{1-3z},
\end{equation}
where $\tau(\theta)$ is the unique real number for which          
$\rho(f_\theta)=\theta$.
$f_\theta$ has a degree three critical point at $1$, and we denote by $W$
the component of $f_\theta^{-1}({\Bbb D})$ not contained in $\Bbb D$. 
For a point $\zeta\in {\Bbb C}$, $f^i(\zeta)=1$ we denote by $W(\zeta)$
the $f^i$-preimage of $W$ attached to $\zeta$, and we call $\zeta$
the {\it root point of $W(\zeta)$}.

The following theorem establishes the connection between $f_\theta$ and 
the quadratic polynomial $P_\theta$.
\begin{thm}[Douady, Ghys, Herman, Shishikura]
\label{surg}
Let $\theta$ be of bounded  type, and let $\Delta$ denote the Siegel
disc of $P_\theta$. Then there exists a quasi-conformal homeomorphism
$\phi:\bar{{\Bbb C}}\to\bar{{\Bbb C}}$,
conformal on the immediate basin of infinity, such that
$\phi({\Bbb D})= \Delta$, and 
$\phi\circ f_\theta=P_\theta\circ\phi$ on ${\Bbb C}\setminus {\Bbb D}$.
\end{thm}
Define $J_\theta=J(f_\theta)\setminus (\cup_{n\geq 0}f_\theta^{-n}(W)\cup
{\Bbb D})$.

\begin{figure}[htbp]
\centerline{\psfig{figure=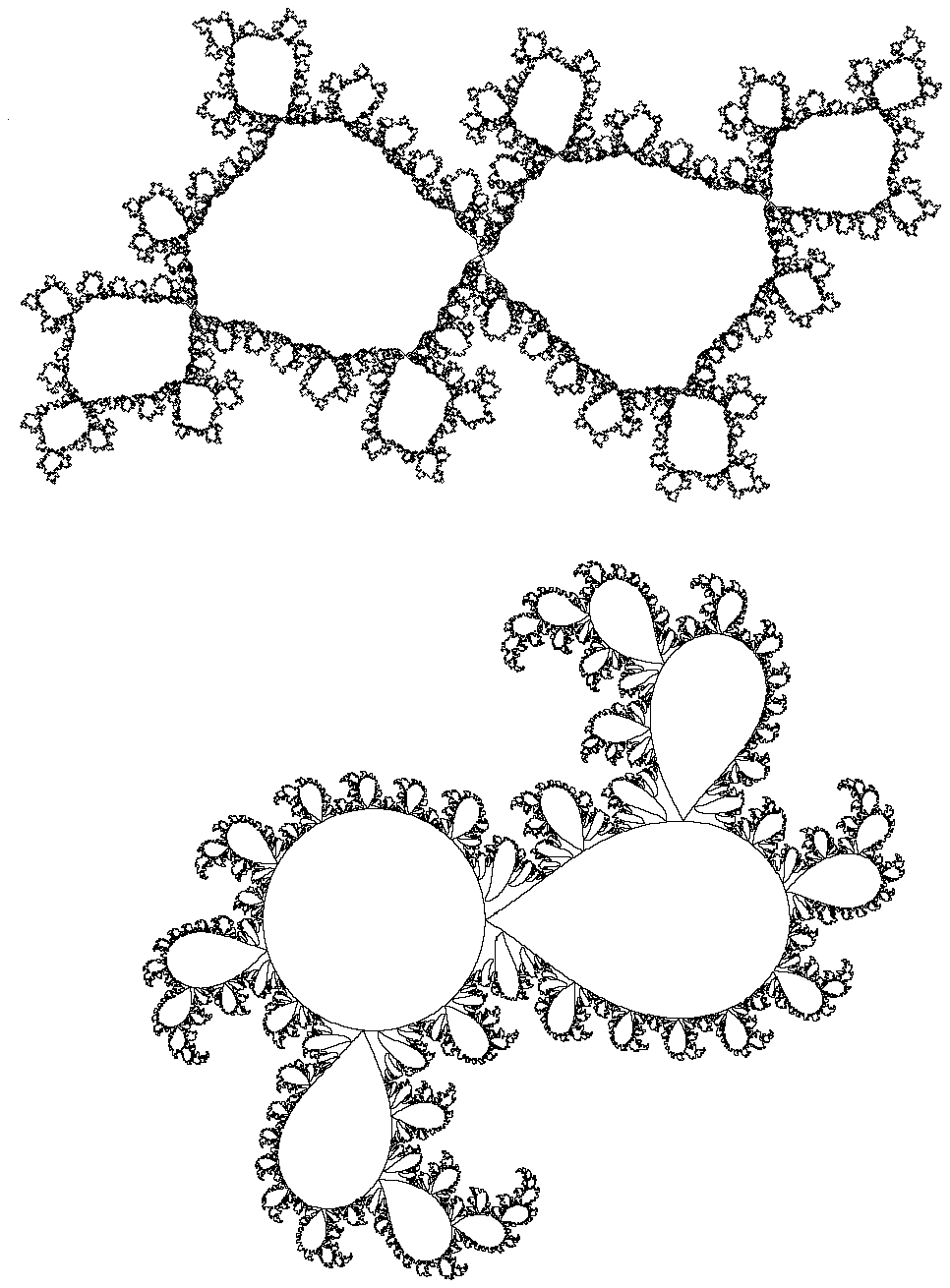,width=13.5cm}}
\vspace{0.7cm}
\caption[surgpic]{A quadratic Siegel Julia set ( above )
and the corresponding set $J_\theta$}
\oplabel{surgpic}
\end{figure}}

Note that for $\phi$ as in \thmref{surg}, one has 
$J(P_\theta)=\phi(J_\theta)$.
Petersen \cite{P} proved the following result:
 
\begin{thm}[Petersen, \cite{P}]

\label{lc1}

Let $P_\theta$ be as above. If $\theta$
is irrational of bounded type, then the Julia set $J(P_\theta)$ is
locally connected and has Lebesgue measure zero.

\end{thm}

He actually gave a proof of another result which together with \thmref{surg}
implies the above theorem.

\begin{thm}[Petersen, \cite{P}]
\label{lc2}
For any irrational rotation number $\theta$  
the set $J_\theta$ is 
locally connected and has Lebesgue measure zero.
\footnote{In \cite{P} Petersen proved the measure zero part of \thmref{lc2}
for $\theta$ of bounded type only. The argument for the general case was
suggested to us by M. Lyubich.}
\end{thm}

In what follows we present a new proof of \thmref{lc2},
 which bears a strong analogy to the proofs of local connectivity of the Julia sets 
for quadratic polynomials (see \cite{HJ,McM2}).  

We give a brief outline of the proof. First, following Petersen
we construct a sequence of  "puzzle-pieces". By the
complex bounds they shrink down to the critical point $1$. This yields the 
basis of connected neighborhoods around the critical point. We then
present a "spreading around" argument to prove local connectivity 
at any other point of the set $J_\theta$. 

\comm{
M. Lyubich has pointed out  that the same argument proves the measure zero 
statement after we notice that each of the above neighborhoods contains a 
definite amount of "empty space", so no point can be a Lebesgue 
 density point of
$J_\theta$.
}

\subsection{Complex bounds}

Let the domain $G={\Bbb C}\setminus \{0,\infty\}$. 
Consider the universal covering $\kappa : {\Bbb C}\to G$,
$\kappa:z\mapsto e^{iz}$. Let $\bar f:{\Bbb C}\to {\Bbb C}$ denote the 
lift of the map 
$f$ to this covering.

The map $\bar f$ belongs to an Epstein class, and the cubic estimate
(\ref{cubic-estimate}) holds for this map. 

Let $D\subset G$ be a domain such that $\operatorname{cl}D\cap T\subset 
[f^{q_2}(1),f^{q_1-q_2}(1)]$. By remark \ref{choice-domain} the domain
$D_0$ in the Proposition \ref{cubic-estimate} for the map $\bar f$
can be chosen so that $\kappa(D_0)\supset D$.

Using the bound on the derivative of the exponential map in the
domain $D_0$,
we obtain the estimate (\ref{cubic-estimate}) for the map $f$ 
itself in the domain $D$, with constants depending only on the choice of
the domain.

\subsection{Construction of "puzzle-pieces"}
We construct a sequence of puzzle-pieces around the
critical point following Petersen.


%
%
\comm{
Let $1\in \gamma_1\subset \partial W$ be on the same side of $W$ as 
$f^{q_n+q_{n+1}}(1)$, such that 
$f^{q_n+q_{n+1}}(\gamma_1)=[1,f^{q_n+q_{n+1}}(1)]$.
Note that $\gamma_1\subset  D_n\subset \Delta_n$.

\realfig{circ4}{circ4.eps}{}{9cm}

Let $\gamma_2\cap\gamma_1\neq \emptyset$ and $f^{q_n+q_{n+1}}(\gamma_2)=
\gamma_1$. Inductively define $\gamma_i\cap\gamma_{i-1}\neq \emptyset$,
and $f^{q_n+q_{n+1}}(\gamma_i)=\gamma_{i-1}$.

The curves $\{\gamma_i\}$ are all contained in $D_n$, and they converge
to a point $\beta_n$ which is fixed by the map $f^{q_n+q_{n+1}}$.
Let $\displaystyle\gamma=\cup_{i\geq 1}\gamma_i\cup \beta_n$.
Let $\Gamma\equiv \Gamma_n=f^{-q_n}(\gamma)\cup[f^{-q_n}(1),1]\subset V_n$, 
and 
let $\Gamma'\equiv \Gamma'_{n+1}\subset V_n$ be the symmetric
of $\Gamma_{n+1}$, i.e. $f_\theta(\Gamma_{n+1})=
f_\theta(\Gamma')$.
Note that both $\Gamma$ and $\Gamma'$ are contained in $\Delta_n$.

Denote the endpoints of $\Gamma$ and $\Gamma'$ by $\nu$ and $\nu'$
respectively. Let $R$ and $R'$ be the external rays landing on $\nu$ and 
$\nu'$ such that any other external ray landing on one of these points 
is contained in the sector $R\cup\Gamma\cup\Gamma'\cup R'$.
Finally let $E\subset \Delta_n$ be a piece of an equipotential
cut out by the rays $R$ and $R'$.

We let the "puzzle-piece" $P_n$ be the closed region cut out by
the curves $\Gamma,\Gamma',R,R',$ and $E$ (cf. \figref{circ4}).
}
%
%

Let $\gamma'_0\subset \partial W$, $\gamma_0\subset \partial W$,
 and $f(\gamma'_0)=[f(1),1]$, $f(\gamma_0)=T\setminus [f(1),1]$.
Let $\gamma'_1\cap \gamma'_0\neq \emptyset$, $ \gamma_1\cap \gamma_0\neq
\emptyset$, and $f(\gamma'_1)=\gamma'_0$, $f(\gamma_1)=f(\gamma_0)$.
Inductively define $\gamma'_i\cap\gamma'_{i-1}\neq \emptyset$,
$\gamma_i\cap\gamma_{i-1}\neq \emptyset$, and 
$f(\gamma'_i)=\gamma'_{i-1}$, $f(\gamma_i)=f(\gamma_{i-1})$.

The curves $\gamma'_i$, $\gamma_i$ converge to a repelling fixed point $\beta$
of the map $f$.

Let $ \Gamma'=\cup_i \gamma'_i \cup \beta$,
 $\Gamma= \cup_i \gamma_i \cup \beta$, and let $\hat\Gamma=f^{-1}(\Gamma')$,
$\hat\Gamma\ni f^{-1}(1)$. Denote by $R$ the external ray of external argument
$0$ landing at $\beta$ and let $R'$ be its preimage landing at the end point
of $\hat\Gamma$. Finally, let $E$ be an equipotential. 

We let the puzzle-piece $P_0\supset W$ be the closed domain cut out by the 
curves $R\cup\Gamma\cup[1,f^{-1}(1)]\cup\hat\Gamma\cup R'$ and $E$ (cf.
\figref{puzzle}).

\realfig{puzzle}{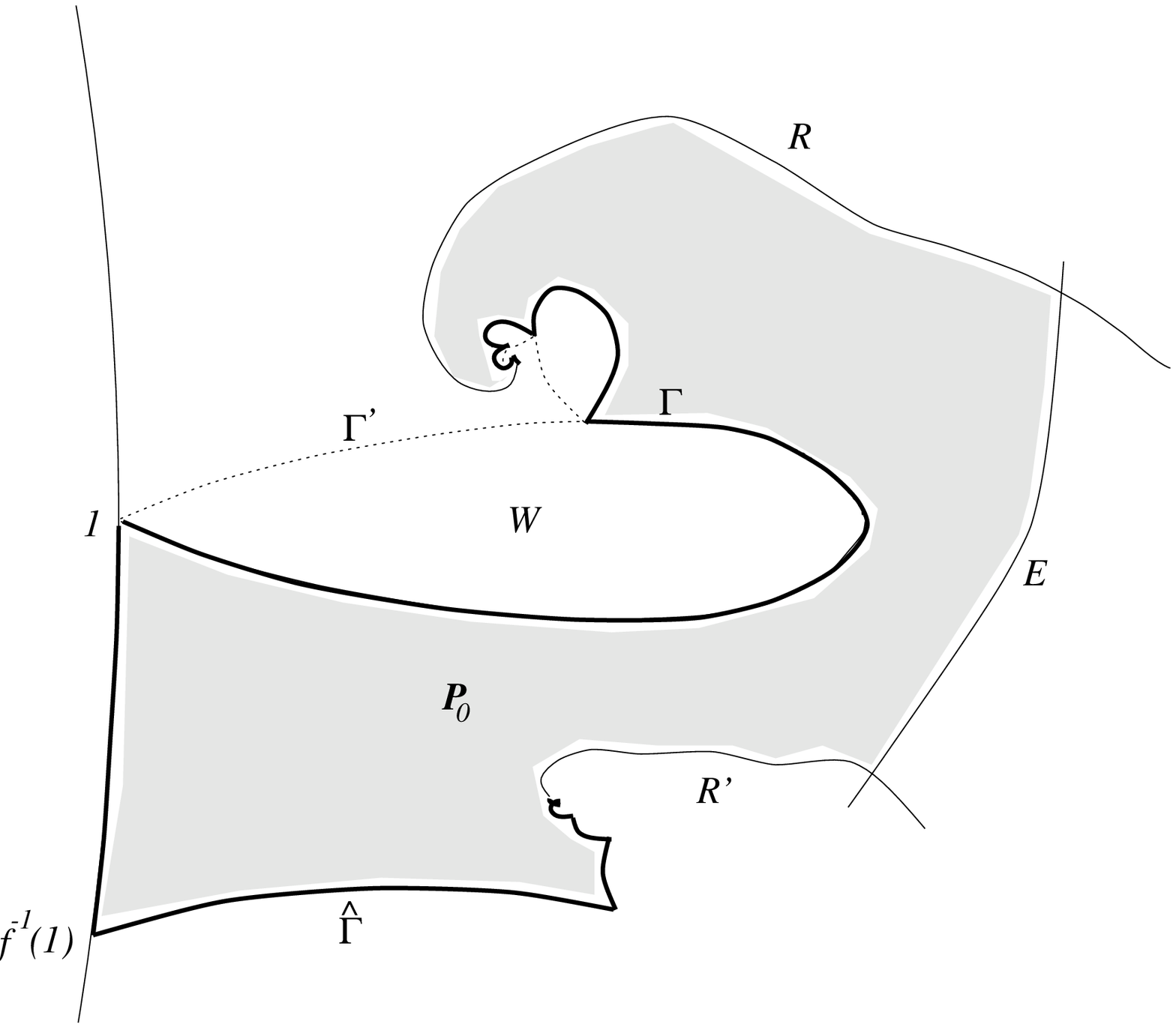}{}{10cm}

Let $P_1$ be the pullback of the domain $P_0$ corresponding to the 
inverse orbit $[1,f^{q_1}(1)]\subset P_0, [f^{-1}(1),f^{q_1-1}(1)],\ldots,
[f^{-q_1}(1),1]$, and inductively, let $P_n$ be the pullback of
$P_{n-1}$ corresponding to the orbit $[1,f^{q_n}(1)]\subset P_{n-1},
[f^{-1}(1),f^{q_n-1}(1)],\ldots,[f^{-q_n}(1),1]$.

Note that $P_n\cap T=[f^{-q_n}(1),1]$, $f^{q_n}(P_n\cap \partial W)=
[f^{q_n}(1),f^{-q_{n-1}}(1)]$, and $f^{q_n+q_{n-1}+q_{n-2}}(P_n\cap \partial W 
(f^{-q_n}(1))=[1,f^{q_{n-1}+q_{n-2}}(1)]$.

By construction we immediately have the following
\begin{prop}
\label{pp1}

The intersection $P_n\cap J_\theta$ is connected.
\end{prop}

By \'Swiatek-Herman real a priori bounds the intervals $[f^{-q_n}(1),1]$ and
 $[1, f^{-q_{n-1}}(1)]$ are 
$K-$commensurable, with a universal constant $K$ for 
sufficiently large $n$.

By the cubic  estimate \ref{cubic-estimate}, 
\begin{center}
$\displaystyle
\diam P_n\leq (C_1\sqrt[3]{\frac{\diam P_{n-1}}{|[f^{-q_{n}}(1),1]|}}+C_2)
\cdot |[ f^{-q_{n}}(1),1]|.$
\end{center}
Hence, if $\displaystyle \frac {\diam P_{n-1}}{|[1, f^{-q_{n-1}}(1)]|}>K_1$
 for a large
 $K_1$, then $\displaystyle\frac {\diam P_{n}}{|[f^{-q_{n}}(1),1]|}<\frac{1}{2
}\cdot
\frac {\diam P_{n-1}}{|[1, f^{-q_{n-1}}(1)]|}$. 
It follows that for all sufficiently large $n$, the piece $P_n$ is
$K_2-$commensurable with $[f^{-q_n}(1),1]$ with a universal constant $K_2$.

Together with Proposition \ref{pp1} this implies
\begin{prop}
\label{lcat1}
The set $J_\theta$ is locally connected at the critical point $1$.
\end{prop}

\begin{lem}
\label{diam2}
$P_n$ contains a Euclidean disc $B$ with $\diam B>K \diam P_n$,
  for a universal constant $K>0$.
\end{lem}
\begin{pf}
Note that by construction, $W(f^{-q_{n+2}}(1))\subset P_n$. The claim 
now easily follows.
\end{pf}

\subsection{``Spreading around'' argument}

Choose any $z\in J_\theta$.

Assume first that there exists $n$ such that $f^{i}(z)\notin P_k$ for
any $i\geq 0$ and $k\geq n$.
As $f$ has an irrational rotation number on the circle this implies that
the forward orbit $z_0\equiv z,z_1\equiv f(z), z_2\equiv f(z_1),\ldots$
does not accumulate on the circle, i.e. there exists $\eps>0$, such that $z_i\in
V_\eps\equiv \{\zeta, \; |\zeta|>1+\eps\}$.

As the set $J_\theta$ is locally connected at the critical point $1$, there
exist two external rays $r_1$ and $r_2$ landing at this point on different
sides of $W$.

For  a point $\zeta$ in the inverse orbit of the critical point 
let $r_1(\zeta)$, $r_2(\zeta)$ be the preimages of the 
rays $r_1$ and $r_2$, landing on $\zeta$. Let $L_\zeta$ be the ``limb''
of the set $J_\theta$, cut  out by the rays $r_1(\zeta), r_2(\zeta)$.

Let $a$ be an accumulation point of the sequence $\{ z_k\}$. Choose a limb 
$L\equiv L_\zeta$ containing $a$, with $L\cap T=\emptyset$. Denote by $k_n$ the
moments when $z_{k_n}\in L$. Let $L_n\ni z$ be the pullback of $L$
corresponding to the backward orbit $z_{k_n},z_{k_n-1},\ldots,
z_{1},z_0\equiv z$. We use the following general lemma to assert that
$\diam (L_n)\to 0$.

\begin{lem}[\cite{L3}, Prop. 1.10]
Let $f$ be a rational  map. Let $\{f^{-m}_i\}$ be a family of 
univalent branches of the inverse functions in a domain $U$.
If $U\cap J(f)\neq \emptyset$, then for any $V$ such that $\cl V\subset U$, 
\begin{center}
$\diam(f^{-m}_i V)\to 0$.
\end{center}
\end{lem}

This yields the desired nest of connected neighborhoods around $z$,
and we are done.

Now let $z_k$ be the first point in the orbit
$z_0,z_1, z_2,\ldots$ contained in the
piece $P_n$.
Denote by 
\begin{equation}
\label{P-orbit} 
\Pi_0\equiv P_n,\Pi_{-1},\ldots,\Pi_{-k}
\end{equation} 
the preimages of $P_n$ corresponding to the inverse orbit $z_k,z_{k-1},\ldots,
z_0$.

\begin{lem}
\label{once}

The inverse orbit (\ref{P-orbit}) hits the critical point $1$ at most once.
\end{lem}
\begin{pf}
To be definite, assume that $P_n$ is above the critical point $1$.
Note that if $\Pi_{-i}\cap T=\emptyset$ for some $i\leq q^{n+1}$, then
the inverse orbit (\ref{P-orbit}) never hits the critical point.
Otherwise, denote by $A$ and $B\equiv P_{n+1}$ the ``above'' and ``below''
$f^{q_{n+1}}$-preimages of $P_n$, i.e. $f^{q_{n+1}}(A)=P_n$, 
$A\cap T\neq\emptyset$, and $A$ is above $1$, and similarly for $B$.
Notice  that $A\cap T=[f^{-q_{n+1}-q_{n}}(1),1]{\subset}
[f^{-q_{n}}(1),1]$. Let $L_1=P_n\cap \partial W$, and $L_2=A\cap W$, then
\begin{eqnarray*}
f^{q_{n-1}+q_n+q_{n+1}}(L_1)&=& f^{q_{n+1}}([f^{q_{n-1}+q_n}(1),1])
=[f^{q_{n-1}+q_n+q_{n+1}}(1),f^{q_{n+1}}(1)]\\
&\supset&
[f^{q_{n-1}+q_n}(1),f^{q_{n-1}+q_n+q_{n+1}}(1)]
=f^{q_{n-1}+q_n}([1,f^{q_{n+1}}(1)])\\
&=&f^{q_{n-1}+q_n+q_{n+1}}(L_2).\\
\end{eqnarray*}
 Thus $L_1\supset L_2$,
and as two different preimages of $W$ cannot cross, it follows that $A\subset
P_n$.
Hence $\Pi_{-q_{n+1}}\neq A$.

Finally, denote by $B'\cap T\neq \emptyset$ the $f^{q_n}$ preimage of $B$,
$f^{q_n}(B')=B$.
$B'\cap T=[f^{-q_n}(1),f^{-q_n-q_{n+1}}(1)]\subset[f^{-q_n}(1),1]$. 
Let $L_2=P_n\cap W(f^{-q_n}(1))$, and $L_3=B'\cap W(f^{-q_n}(1))$,
then
\begin{eqnarray*}
f^{q_n+{q_{n-1}}+q_{n-2}}(L_2)&=&[1,f^{q_{n-1}+q_{n-2}}(1)]\\
&\supset& [f^{q_{n-1}+q_{n-2}-q_n-q_{n+1}}(1),f^{q_{n-1}+q_{n-2}}(1)]\\
&=&f^{q_n+{q_{n-1}}+q_{n-2}}(L_3).
\end{eqnarray*}
Therefore, $L_2\supset L_3$, and 
$B'\subset P_n$.
\\
Thus, $\Pi_{-q_{n+1}-q_{n}}\cap T=\emptyset$ and the claim follows.
\end{pf}

Now let $-k\leq -i\leq 0$ be the first moment such that
 $\Pi_{-i}\cap T=\emptyset$. By real bounds there exists an annulus $A$ around
$\Pi_{-(i-1)}$ with $\mod A>K$ for a universal $K$, such that the critical value
$f(1)$ is outside of $A$. Therefore, there exists an annulus with modulus $K$
around $\Pi_{-i}$, such that the whole postcritical set of $f$ is outside of
it.

By \lemref{once} and Koebe theorem, the inverse branch 
$f^{-k}:P_n\to \Pi_{-k}\equiv P_n(z)\ni z$ has bounded distortion
 on the piece $P_n$.
As $P_{n}(z)$ cannot contain a disc of definite radius, and by \lemref{diam2},
$\diam(P_{n}(z))\to 0$, which yields the desired nest of connected 
neighborhoods around $z$.

\subsection{Proof of measure zero statement}

The following proof was suggested by M.~ Lyubich. 
First consider the set of points $J_1\subset J_\theta$, $J_1\equiv
\{ \zeta\in J_\theta |\; \exists n, \; f^i(\zeta)\notin P_k,\;
\forall i\geq 0, k\geq n\}$.
The set $J_1$ has zero Lebesgue measure
by a  theorem of  M.Lyubich
(\cite{L1}).

For a  piece $P_n$ denote by $P'_n$ the symmetric piece with respect to
the circle $T$. Let $Q_n\equiv P_n\cup P'_n$. Note, that $\operatorname{int}
P'_n\cap
J_\theta =\emptyset$.
Consider now a point $z\in J_\theta\setminus J_1$. Let $k$ be the first moment
when the orbit $z_0\equiv z, z_1,z_2,\ldots$ enters the piece $P_n$.
By the same argument as above, there is a piece $Q_n(z)$ around $z$, which
is a bounded distortion pull-back of $Q_n$. The diameters of sets $Q_n(z)$
tend to zero.

Thus we obtain a sequence of balls $B_n\supset Q_n(z)$ around $z$,
such that by \lemref{diam2},
\begin{center}
$\displaystyle\frac{{\meas}(B_n\setminus J_\theta)}
{{\meas}(B_n)}>\delta>0$.
\end{center}
It follows that $z$ is not a point of Lebesgue density of set $J_\theta$, and
this completes the proof of \thmref{lc2}.

\end{document}